\documentclass[12pt,letterpaper]{amsart}

\usepackage{amssymb,amscd}

\newcommand{\p}{{\mathcal P}}

\newcommand{\oo}{{\mathcal O}}

\newcommand{\F}{{\mathbf F}}

\newcommand{\Q}{{\mathbf Q}}

\newcommand{\spec}{{\text{Spec}\,}}

\newcommand{\vsp}{\vspace{8pt}}

\newcommand{\vtsp}{\vspace{16pt}}
\newcommand{\rarr}{{\rightarrow}}

\newcommand{\mymod}[2]{{ #1 \: (\bmod \: {#2})}}

\newtheorem{thm}{Theorem}
\newtheorem{lem}{Lemma}

\newtheorem*{corr}{Corollary}

\theoremstyle{remark}
\newtheorem{remark}{Remark}

\newtheorem*{ack}{Acknowledgment}

\newcommand\gal{{ \mbox{\rm Gal}}}

\newcommand{\ov}[1]{{\overline{{#1}}}}

\usepackage{fullpage}

\newcount\timehh\newcount\timemm
\timehh=\time 
\divide\timehh by 60 \timemm=\time
\begin{document}

\title{Quadratic twists of pairs of elliptic curves}

\author{Siman Wong
  {
    \protect \protect\sc\today\ -- 
    \ifnum\timehh<10 0\fi\number\timehh\,:\,\ifnum\timemm<10 0\fi\number\timemm
    \protect \, \, \protect \bf DRAFT
  }
}

\address{Department of Mathematics \& Statistics, University of Massachusetts.
	Amherst, MA 01003-9305 USA}

\email{siman@math.umass.edu}





\begin{abstract}
Given two elliptic curves defined over a number field $K$, not both with
$j$-invariant zero,
we show that there are infinitely many $D\in K^\times$ with pairwise distinct
image in
$
K^\times/{K^\times}^2
$,
such that the quadratic twist of both curves by $D$ have
positive Mordell-Weil rank.  The proof depends on relating the 
values of pairs of cubic polynomials to rational
points on another elliptic curve, and on a fiber product construction.
\end{abstract}

\maketitle


\section{Introduction}

Fix two elliptic curves over $\Q$ with coprime conductors.  Then the
parity conjecture predicts that there are infinitely many square-free
integers $D$ so that the quadratic twist of both curves by $D$ have positive
Mordell-Weil rank.  This argument no longer applies if the curves have
opposite root numbers and if their conductors have the same square-free part,
not to mention the fact that it is based on a deep conjecture.
Furthermore, Rohrlich informed us that there exist number fields $K$ and
elliptic curves $E$ over $K$
for which every quadratic twist of $E$ over $K$ has even analytic rank
\cite{rohrlich}.
%
%
%
In this paper we give an unconditional construction of
simultaneous positive rank twists under a mild restriction on $j$-invariants.

\vsp

\begin{thm}
    \label{thm:pair}
For any pair of elliptic curves defined over a number field $K$, not both
with $j$-invariant zero, there exist infinitely many $D\in K^\times$
with pairwise distinct image in
$
K^\times/{K^\times}^2
$,
such that the quadratic twist by $D$ of both curves have positive Mordell-Weil
rank over $K$.
\end{thm}

\vsp

\begin{remark}
The elements $D$ in the theorem are values of a cubic polynomial at
the $x$-coordinate of multiples of a non-torsion point on
another elliptic curve.  In particular, they form a thin
set, while the parity conjecture when applicable yields a set of positive
density.
\end{remark}

\vsp

Apply theorem \ref{thm:pair} to the case where one curve in the pair is the
quadratic twist of the other one and we deduce the following result.

\vsp

\begin{corr}
    \label{corr:simu}
Let $E$ be an elliptic curve defined over a number field $K$ with non-zero
$j$-invariant.  Then for any 
$
\delta\in K^\times
$,
there exist infinitely many $D\in K^\times$ with pairwise distinct image in
$
K^\times/{K^\times}^2
$,
such that the quadratic twist
of $E$ by both $D$ and $D\delta$ have positive Mordell-Weil rank over $K$.
      \qed
\end{corr}

\vsp

\begin{remark}
In general we do not know that the elements $D$ furnished by the corollary
are coprime to $\delta$.  This extra requirement does hold trivially when
$\delta$ is irreducible.
\end{remark}

\vsp

In the case of $j$-invariant zero we have the following partial result.

\vsp

\begin{thm}
    \label{thm:zero}
Let $E_1, E_2$ be elliptic curves defined over a number field $K$ with
$j$-invariant zero.  Then there exists
$
\lambda\in K^\times
$
such that, if we denote
by $E_i'$ the sextic twist of $E_i$ by $\lambda$, then there are infinitely
many $D\in K^\times$ with pairwise distinct image in
$
K^\times/{K^\times}^2
$,
such that the quadratic twist by $D$ of both
$E_1'$ and $E_2'$ have positive Mordell-Weil rank over $K$.
\end{thm}

\vsp

At most finitely many quadratic twists of a fixed elliptic curve over $K$ has 
$K$-rational torsion
points of order $>2$.  To find positive rank twists of 
$
y^2 = f(x)
$
it then suffices to find $x_0\in K$ so that
$
f(x_0)\in K - K^2
$.
Given another elliptic curve
$
y^2 = g(x)
$,
to find simultaneous positive rank twists for both curves we are then led to
 seek conditions
on $f$ and $g$ so that
\begin{itemize}
\item[(i)]
the cubic equation
$
f(x) = g(y)
$
defines another elliptic curve $E'/K$;
\item[(ii)]
we can construct a point $P'$ of infinite order on $E'(K)$;
and
\item[(iii)]
we evaluate
$
f(x) = g(y)
$
at multiples of $P'$ to generate infinitely many values with pairwise distinct
image in
$
K^\times/{K^\times}^2
$.
\end{itemize}
The solutions to
$
f(x) = g(y)
$
for which this common value has a given image 
$
\lambda\in K^\times/{K^\times}^2
$
turn out to be parameterized by a certain fiber product
$
C_\lambda
$
of elliptic curves; Step (iii) then comes down to showing that every
$C_\lambda$ has geometric genus $>1$.  As for Step (ii),
the cubic equation $f(x) = g(y)$ contains several natural 
solutions  if at least one of
$
y^2 = f(x)
$
and
$
y^2 = g(x)
$
has non-zero $j$-invariant, but we need to show that these points give
rise to non-torsion points on
$
E'(K)
$.
We do that by adjusting the Weierstrass models of the curves.

\vsp

We expect theorem \ref{thm:pair} to hold with no restriction on
$j$-invariant.
Removing this
condition, however, seems non-trivial, cf.~Remark \ref{rem:triple}.
Similarly, theorem \ref{thm:pair} should hold for any finite collection
of elliptic curves, but our argument does not generalize to this setup.

\vsp

\section{A fiber product}
    \label{sec:simu}

Given a pair $(a, b)\in K^2$ with $4a^3 + 27b^2\not=0$, let
$
f_{a, b}(x) = x^3 + ax + b
$,
and denote by
$
E_{a,b}
$
the elliptic curve over $K$ defined by the Weierstrass equation
$
y^2 = f_{a, b}(x)
$.
Given two such pairs $(a, b), (c,d)$,
  denote
  by
$
E' = E'_{a, b, c, d}
$
the projective plane curve over $K$ defined by
\begin{equation}
z^3 f_{a, b}(x/z) = z^3 f_{c, d}(y/z).
          \label{plane}
\end{equation}
It contains the rational point
$
[x: y: z] = [1: 1: 0]
$
and, provided that $a\not=c$, another point
$
P' = P'_{a, b, c, d} = [b-d: b-d: c-a]
$.
In the next section we will study conditions under which $E'$ is non-singular
and for which
$
P'
$
is a non-torsion point.  For the rest of this section we will investigate
Step (iii), but we should point out that
$
P'
$
becomes  $[1: 1: 0]$ when
$
c-a=0
$.
This happens notably when $a=c=0$, i.e.~when both curves have $j$-invariant
zero, and this turns out to be the source of the restriction on
$j$-invariant in the
theorem.  If we try to generalize (\ref{plane}) by using 
$
z^3 f_{a, b}(x/z) = \mu^2 z^3 f_{c, d}(y/z)
$,
the new equation seems to have no natural rational point unless
$
\mu
$
is a perfect cube, in which case we are essentially back to (\ref{plane}).

\vsp

Denote by
$
E_\lambda
$
the twisted elliptic curve
$
\lambda y^2 = f_{a, b}(x)
$,
and by
$
\varphi_\lambda: E_\lambda \rarr\mathbf P^1
$
the projection defined by the affine map
$
(x, y) \mapsto x
$.
Denote by
$
\psi: E' \rarr \mathbf P^1
$
the projection given on the affine curve $f_{a, b}(x) = f_{c, d}(y)$ by
$
(x, y) \mapsto x
$.
Denote by $C_\lambda/K$ the fiber product
$
E_\lambda \times_{\mathbf P^1} E'
$
defined via $\varphi_\lambda$ and $\psi$.  This is a projective curve over $K$, and
it has infinitely many $K$-rational points precisely if there are infinitely
many pairs
$
(x_0, y_0)\in K^2
$
such that the common value
$
f(x_0) = g(y_0)
$
is $\lambda$ times a perfect square in $K$.  To handle Step (iii) it then
suffices to show that every $C_\lambda$ has geometric genus $>1$
(we will clarify this in the proof of theorem \ref{thm:pair}).

\vsp

\begin{lem}
      \label{lem:genus}
Let
$
(a, b)
$
and
$
(c,d)
$
be as above.  If $E'$ is non-singular and if $\lambda\not=0$,
then the
fiber product $C_\lambda$ is non-singular and has genus $>1$.
\end{lem}

\begin{proof}
The branched locus of $\varphi_\lambda$ consists precisely of the roots of
$
f_{a, b}(x)
$
(where as usual we view $\ov{K}$ as
$
\mathbf A^1_\ov{K}\subset\mathbf P^1_\ov{K}
$)
plus the point $\infty\in\mathbf P^1$.  The fiber of $\psi$ above $\infty$
consists of those solutions to (\ref{plane}) with $z=0$; there are three 
distinct ones, so the triple cover $\psi$ is not branched over
$
\infty
$.
Now, $\psi$ is ramified above a finite point
$
x_0\in\ov{K}
$
precisely when
$
f_{c, d}(y) - f_{a, b}(x_0)
$,
as a cubic polynomial in $y$, is not separable.  This is equivalent to
$
4c^3 + 27(d - f_{a, b}(x_0))^2 = 0
$.
Since $f_{c, d}$ is separable, it means that
$
f_{a, b}(x_0)\not=0
$.
Thus the branched loci of $\varphi_\lambda$ and $\psi$ are disjoint.
Furthermore,  $E'$ and $E_{a, b}$ are both non-singular, so the fiber product
$C_\lambda$ is also non-singular.   To compute its genus we
can then apply the Riemann-Hurwitz formula.

\vsp

Denote by
$
\pi:C_\lambda\rarr\mathbf P^1
$
the projection coming from the fiber product.  We saw that  the double cover
$\varphi_\lambda$ is
branched above four points
in $\mathbf P^1$ (since $f_{a, b}$ is separable), and that the triple cover
$\psi$
is unramified above  each of these points.  This gives $12$ points on 
$
C_\lambda(\ov{K})
$
at which $\pi$ has ramification index $2$.  By the Riemann-Hurwitz formula,
the ramification of $\pi$ above the branched points of $\varphi_\lambda$ alone implies
that $C_\lambda$ has genus $\ge 1$.  But
 $\psi$ and $\varphi_\lambda$
have disjoint branched loci, so $\pi$ has additional ramification,
whence $C_\lambda$ must have genus $>1$.
\end{proof}

\vsp

\section{Non-torsion points}
    \label{sec:tor}

We begin with a criterion for $E'_{a, b, c, d}$ to be non-singular.
Using the command {\tt Weierstrassform} in the computer algebra system
{\tt MAPLE}, we find that the transformation
\begin{eqnarray}
X
&=&
3x^2+a+3yx+3{y}^{2}+c
     \label{xx}
\\
Y
&=&
-3ya-6ax-3cx- 9b/2+3cy+9d/2-9y{x}^{2}-9{y}^{2}x-9{x}^{3}
     \nonumber
\end{eqnarray}
takes $f_{a, b}(x) = f_{c, d}(y)$ to
\begin{equation}
E'':
Y^2 = X^{3}-3acX
-
{a}^{3} - {c}^{3} - 27(b-d)^2/4
        \label{ee}
\end{equation}
(actually (\ref{xx}) is the negative of that furnished by {\tt MAPLE}; we make
this adjustment so that $E''$ takes on the usual Weierstrass form).
Under this transformation, the point
$
P'_{a, b, c, d}
$
becomes
\begin{equation}
P''
=
P''_{a, b, c, d}
=
\Bigl(
  \frac{9(b-d)^2  + (a-c)^2(a+c)}
       {(a-c)^{2}},
\:\:
\frac{9(b-d)(6(b-d)^2 + (a-c)^2(a+c))}
     {2(a-c)^{3}}
\Bigr)
               \label{point}
\end{equation}
(we will address the vanishing of $a-c$ later on).  So if the discriminant
of the cubic on the right side of (\ref{ee}), namely
\begin{eqnarray}
108 a^3 c^3
-27( 4a^3 + 4c^3  + 27(b-d)^2 )^2 / 16,
        \label{disc}
\end{eqnarray}
is non-zero, then $E''$ is an elliptic curve, and hence so does $E'$.
Note that these
{\tt MAPLE} computations are purely symbolic and hence applies to all
sufficiently large characteristics.

\vsp

For any number field $k$, denote by $\oo_k$ the ring of integers of $k$,
and by
$
\F_{\mathfrak p}
$
the residue field of
$
\mathfrak p\in\text{Spec}\, \oo_k
$.
For any $u\in K^\times$, we write $\mathfrak p\nmid u$ if $u$ is a
$\mathfrak p$-adic unit.

\vsp

\begin{lem}
    \label{lem:tor}
Given two elliptic curves over $K$
which are not isomorphic over $K$ and not both with $j$-invariant zero, 
we can find Weierstrass equations
$
E_{a, b}
$
and
$
E_{c, d}
$
for them so that $E''_{a, b, c, d}$ is an elliptic curve, and that
$
P''_{a, b, c, d}
$
is a non-torsion point in $E''_{a, b, c, d}(K)$.
\end{lem}

\begin{proof}
Fix Weierstrass equations
$
E_{a, b}
$
and
$
E_{c, d}
$
for these two curves.  Without loss of generality, suppose 
$
E_{a, b}
$
has non-zero $j$-invariant; equivalently, $a\not=0$.  

\vsp

First, suppose $b=0$.  Fix $\mathfrak p\in\text{Spec}\, \oo_K$ with
$
\mathfrak p \nmid 6a
$,
and fix a non-zero element $\pi\in\mathfrak p$.  Set
$
\ov{c} = \pi^4 c
$
and
$
\ov{d} = \pi^6 d
$.
Then
$
E_{\ov{c}, \ov{d}}
$
defines the same elliptic curve over $K$ as $E_{c, d}$, and
\begin{itemize}
\item
the discriminant of 
$
E''_{a, b, \ov{c}, \ov{d}}
$
is non-zero modulo $\mathfrak p$, so
$
E''_{a, b, \ov{c}, \ov{d}}
$
is non-singular over $\F_{\mathfrak p}$, and hence over $K$;
\item
$
a\not\equiv\mymod{\ov{c}}{\mathfrak p}
$,
so
$
P''_{a, b, \ov{c}, \ov{d}}
$
is well-defined over $\F_{\mathfrak p}$, and hence over $K$; and
\item
$
P''_{a, b, \ov{c}, \ov{d}}
$
is a $2$-torsion point in
$
E''_{a, b, \ov{c}, \ov{d}}(\F_{\mathfrak p})
$
but not in
$
E''_{a, b, \ov{c}, \ov{d}}(K)
$.
\end{itemize}
Apply Merel's theorem on torsion points \cite{merel} and we see that, 
providing that the residual characteristic of $\mathfrak p$ is large enough
(depending on $K$),
$
P''_{a, b, \ov{c}, \ov{d}}
$
is a non-torsion point in 
$
E''_{a, b, \ov{c}, \ov{d}}(K)
$.
A similar argument covers the case $d=0$.  For the rest of the proof we will
take $bd\not=0$.

\vsp

\begin{lem}
    \label{lem:claim}
Suppose $bd\not=0$.  Then lemma \ref{lem:tor} would follow if there exists a
$
\lambda\in K^\times
$
and
$
\mathfrak p\in\spec \oo_K
$
not dividing  $6, a, b, d$ and $4a^3+27b^2$, such that
$
b\equiv\mymod{\lambda^6 d}{\mathfrak p}
$
and
$
a\not\equiv\mymod{\lambda^4 c}{\mathfrak p}
$
{\rm(}in which case $\mathfrak p\nmid \lambda$ as well\/{\rm)}.
\end{lem}

\begin{proof}[Proof of Lemma  \ref{lem:claim}]
With  $\mathfrak p$ and $\lambda$ as above, set
$
\ov{c} = \lambda^4 c
$
and
$
\ov{d} = \lambda^6 d
$.
The hypothesis
$
\mathfrak p\nmid 6(4a^3 + 27b^2)
$
means that the discriminant of
$
E''_{a, b, \ov{c}, \ov{d}}
$
is non-zero modulo $\mathfrak p$, so 
$
E''_{a, b, \ov{c}, \ov{d}}
$
is non-singular over $\F_{\mathfrak p}$, and hence over $K$.  This conclusion remains
true if we replace $\lambda$ by
$
\lambda + \pi
$
for any $\pi\in\mathfrak p$.  So we can choose $\pi$ so that
\begin{equation}
\text{$a\not= \ov{c}$, \:\:
$
b\not= \ov{d}
$,
\:\:
and
\:\:
$6(b-\ov{d})^2 + (a-\ov{c})^2(a+\ov{c}) \not=0$.}
      \label{abc}
\end{equation}
The first two conditions above imply that
$
P''_{a, b, \ov{c}, \ov{d}}
$
is well-defined over $K$, and the last condition says that it is not a $2$-torsion
point in
$
E''_{a, b, \ov{c}, \ov{d}}(K)
$.
Since
$
b\equiv\mymod{\ov{d}}{\mathfrak p}
$
and $\mathfrak p\nmid 2a$, from (\ref{point}) we see that this point reduces to a point
of order $2$ in
$
E''_{a, b, \ov{c}, \ov{d}}(\F_{\mathfrak p})
$.
Apply Merel's theorem as before we see that
$
P''_{a, b, \ov{c}, \ov{d}}
$
is not a torsion point in
$
E''_{a, b, \ov{c}, \ov{d}}(K)
$.
\end{proof}

\vsp

Set
$
\alpha = c/a
$
and
$
\beta = d/b
$;
the congruence condition in lemma \ref{lem:claim} can then be written as
\begin{equation}
\beta \equiv\mymod{\lambda^6}{\mathfrak p}
\:
\text{ and }
\:
\alpha \not\equiv\mymod{\lambda^4}{\mathfrak p}.
      \label{claim}
\end{equation}
Suppose there exist $\lambda, \mu\in K$ such that
$
\beta = \lambda^6
$
and
$
\alpha = \mu^4
$.
Since $E_{a, b}$ and $E_{c, d}$ are not $K$-isomorphic, \cite[p.~49]{joe1}
implies that
$
\lambda^4 \not=\mu^4
$,
whence (\ref{claim}) is satisfied by any
$
\mathfrak p\in\spec\oo_K
$
with
$
\mathfrak p \nmid (\lambda^4 - \mu^4)
$.
From now on we will therefore assume that at least one of
$
K_\alpha =  K(\alpha^{1/4})
$
or
$
K_\beta =  K(\beta^{1/6})
$
is a non-trivial extension of $K$.  
To finish the proof of the lemma, we consider three cases based on conditions
on $K_\alpha$ and
$
K_\beta
$.
In two cases the condition (\ref{claim}) will be satisfied so lemma
\ref{lem:claim}
is
applicable\footnote{in
  what follows, when we can find $\mathfrak p$ satisfying
  (\ref{claim}) we can find infinitely many of them, so the additional
  non-divisibility condition in lemma \ref{lem:claim} is not an issue.
  We will also assume without further
  comment that every $\mathfrak p$  in what follow
  to be unramified in $K_\alpha/K$ and in $K_\beta/K$.};
in
the third case we need to proceed differently.

\vsp

\noindent
\textit{Case I:}
$
K_\alpha \subsetneq K_\beta
$

$
\text{Spec}\, \oo_{K_\beta}
$
has infinitely many $\p$ of degree $1$ over $K$ (i.e.~its
$
K_\beta/K$-norm  is in $\spec\oo_K$); for any such $\p$,
the
congruence
$
x^6 \equiv\mymod{\beta}{\p}
$
is solvable in $\oo_K$.  And if for infinitely many such $\p$, some
$\oo_K$-solution of this congruence is also congruent modulo $\p$ to a
root in $K_\alpha$ of
$
x^4 = \alpha 
$,
say $\alpha_1\in K_\alpha$, then $\alpha_1$ would be an actual sixth
root of $\beta$, contradicting
$
K(\alpha_1 ) \subset K_\alpha \subsetneq K_\beta =  K(\beta^{1/6}) \subset
K(\alpha_1 )
$.
So for infinitely many 
$
\p\in\spec\oo_{K_\beta}$ of degree $1$ over $K$, we can find
$
\lambda\in\oo_K
$
(depending on $\p$) such that
$
\beta \equiv\mymod{\lambda^6}{\p}
$
and
$
\alpha \not\equiv\mymod{\lambda^4}{\p}
$.
Both sides of each of these two congruences are in $K$, so we can
change the modulus from $\p$ to $\p\cap\oo_K$, which is in
$
\spec \oo_K
$
since $\p$ has degree $1$ over $K$, and (\ref{claim}) follows.

\vsp

\noindent
\textit{Case II:}
$
K_\alpha\not\subset K_\beta
$

If $L_\beta$, the Galois closure of $K_\beta /K$, does not contain $K_\alpha$,
then it does not contain the Galois closure of
$
K_\alpha/K
$
either, in which case we can find infinitely many
$
\mathfrak p\in\spec \oo_K
$
such that
$
\spec \oo_{K_\beta}
$
has a degree $1$ prime lying above $\mathfrak p$, but
$
\spec \oo_{K_\alpha}
$
does not.  Once $\mathfrak p$ is chosen, $\lambda$ as in (\ref{claim}) follows
as above.

\vsp

Now, suppose $K_\alpha\subset L_\beta$.  Since $[K_\beta:  K]$
divides $6$, from
$
K_\alpha\not\subset K_\beta
$
and
$
K_\alpha\subset L_\beta
$
we see that
$
[K_\beta:  K] = 3
$
or $6$.  If
$
[K_\beta:  K] = 3
$,
then
$
L_\beta/K
$
is a dihedral extension of degree $6$; since
$
[K_\alpha:  K]
$
divides $4$,  that means 
$
K_\alpha
$
is the unique quadratic subfield of $L_\beta$.  Let
$
\mathfrak p\in \spec\oo_K
$
be unramifed in $L_\beta$, and that its Frobenius conjugacy class is the
class of order $2$ elements in
$
\gal(L_\beta/K)
$.
Then $\mathfrak p$ is inert in $K_\alpha$, and $\oo_{K_\beta}$ has a maximal
ideal $\p$ of degree $1$ lying above $\mathfrak p$.  That means
$
\alpha\not\equiv\mymod{\lambda^4}{\mathfrak p}
$
for any $\lambda\in\oo_K$, while
$
\beta \equiv\mymod{\lambda^6}{\p}
$
has a solution in $\oo_K$.  As before we can replace $\p$ by $\mathfrak p$,
so the condition (\ref{claim}) is satisfied.

\vsp

Next, suppose $K_\alpha \subset L_\beta$ and
$
[K_\beta:  K] = 6
$.
If $K_\beta/K$ is Galois, then from $[K_\alpha:  K]$ dividing $4$ we see
that $K_\alpha/K$ is quadratic, and the argument in the paragraph above
applies.  Now, suppose $K_\beta/K$ is not Galois, which happens if and only 
if $\sqrt{-3}$ is not in $K_\beta$.  That means
$
K_\beta'
$,
the unique cubic subfield
$
K(\beta^{1/3})
$
of $K_\beta$, is not Galois; denote by
$
L_\beta'/K
$
its Galois closure.  This is a diheral extension of degree $6$, and its
unique quadratic subfield is not $K(\sqrt{\beta})$, otherwise $K_\beta/K$
would be Galois.  Thus
$
L_\beta' \cap  K(\sqrt{\beta}) =  K
$
and of course
$
L_\beta = L_\beta'(\sqrt{\beta})
$,
so
$$
\gal(L_\beta/K)  \simeq  \gal(L_\beta'/K) \times \gal( K(\sqrt{\beta})/K).
$$
Denote by
$
\gamma
$
the conjugacy class of elements of 
$
\gal(L_\beta/K)
$
that projects to the class of order $2$ elements in 
$
\gal(L_\beta'/K)
$
and to the trivial class in
$
\gal(K(\sqrt{\beta})/K)
$.
Then for any
$
\mathfrak p_\gamma \in \spec\oo_K
$
which is unramifed in $L_\beta$ and whose Frobenius conjugacy class in
$
\gal(L_\beta/K)
$
is $\gamma$, there is a maximal ideal in each of 
$
\spec \oo_{K_\beta'}
$
and
$
\spec \oo_{K(\sqrt{\beta})}
$
of degree $1$ over $K$ lying above $\mathfrak p_\gamma$, but $\mathfrak p_\gamma$ does not split
completely in the unique
quartic subfield of $L_\beta/K$.  The first statement means that
$
\beta \equiv\mymod{\lambda^6}{\mathfrak p_\gamma}
$
has a solution in $\oo_K$.  We claim that the second statement means that
$
\spec \oo_{K_\alpha}
$
has no maximal ideal of degree $1$ over $K$ lying above $\mathfrak p_\gamma$,
in which case
$
\alpha\not\equiv\mymod{\lambda^4}{\mathfrak p_\gamma}
$
has no solution in $\oo_K$, and the condition (\ref{claim}) is satisfied.

\vsp

Note that $[K_\alpha: K]$ divides $4$ and $K_\alpha\not\subset K_\beta$, so
$
[K_\alpha: K] = 4
$
or $2$.  If
$
[K_\alpha: K] = 4
$,
then $K_\alpha$ is the unique quartic subfield of $L_\beta$, and hence Galois;
the claim then follows immediately from our earlier observation that
$\mathfrak p_\gamma$ does not splitting completely in 
$
K_\alpha/K
$.
If $[K_\alpha:  K]=2$, from
$
K(\sqrt{\beta})\subset K_\beta
$
and
$
K_\alpha\subset L_\beta
$
we see that
$
K_\alpha \not=  K(\sqrt{\beta})
$,
so $\mathfrak p\in\spec\oo_K$ is inert in $K_\alpha/K$ if and only if its Frobenius
in 
$
\gal(L_\beta'/K)
$
is the unique class of order $2$ elements. 
 Recall the definition of $\mathfrak p_\gamma$ and we are 
done.

\vsp

\noindent
\textit{Case III:}
$
K_\alpha = K_\beta
$

Since $[K_\alpha: K]$ divides $4$, $[K_\beta: K]$ divides $6$, and since at
least one of
$
K_\alpha, K_\beta
$
is a non-trivial extension of $K$, that means
$
K_\alpha = K_\beta
$
is a quadratic extension of $K$.  Consequently,
$
\alpha = \alpha_0^2
$
and
$
\beta = \beta_0^3
$
for some $\alpha_0, \beta_0\in K$, with
$
K(\sqrt{\beta_0}) =  K(\sqrt{\alpha_0})\not= K
$;
the equality means that
$
\beta_0 = \alpha_0 \alpha_1^2
$
for some $\alpha_1\in K$, and the inequality means that $\alpha_0\in K$ is
not a square.  Note that
$
\alpha_1=1
$
corresponds precisely to the case where the two curves are non-trivial
quadratic twists of one another.

\vsp

Set
$
\ov{c} = \lambda^4 c /\alpha_1^4 = a \alpha_0^2 \lambda^4 / \alpha_1^4
$
and
$
\ov{d} = \lambda^6 d  /\alpha_1^6 = b\alpha_0^3 \lambda^6
$.
Then
$
E_{c, d}
$
and
$
E_{\ov{c}, \ov{d}}
$
are isomorphic over $K$, and
$$
\begin{array}{lllllllll}
\ov{d} - b
&=&
b(\alpha_0^3\lambda^6 - 1)
&=&
b(\alpha_0\lambda^2 - 1)( (\alpha_0\lambda^2)^2 + \alpha_0\lambda^2 + 1),
\\
\ov{c} - a
&=&
\displaystyle
a\Bigl(
  \frac{\alpha_0^2\lambda^4}{\alpha_1^4} - 1 
\Bigr)
&=&
\displaystyle
\frac{a(\alpha_0 \lambda^2 - \alpha_1^2)(\alpha_0 \lambda^2 + \alpha_1^2)}
     {\alpha_1^4}.
\raisebox{22pt}{}
\end{array}
$$
Since $\alpha_0\in K$ is not a square,
$
\alpha_0\lambda^2 - 1
$
viewed as a polynomial in $\lambda$ is $K$-irreducible.  Thus
$$
\{ \mathfrak p\in\spec\oo_K:
       \text{ $\nu_{\mathfrak p}(\alpha_0\lambda^2 - 1)$ is positive for
             some
	     $\lambda\in\oo_K$}
\}
$$
is infinite, where
$
\nu_{\mathfrak p}
$
denotes an additive $\mathfrak p$-adic valuation for $\oo_K$.  Choose 
$
\mathfrak p\in S
$
for which
\begin{equation}
a, \alpha_0, \alpha_1, \alpha_0^2 \pm 1, \alpha_1^{4}-1
  \text{ and }
\alpha_1^{12}-1
        \label{forwhich}
\end{equation}
are all $\mathfrak p$-adic units, and pick a $\lambda\in\oo_K$ (depending on
$\mathfrak p$) so that
$
\nu_{\mathfrak p}( \alpha_0 \lambda^2 - 1) > 0
$.
Then
\begin{itemize}
\item[(i)]
$
\ov{c}\not\equiv\mymod{a}{\mathfrak p}
$
since
$
\nu_{\mathfrak p}(\alpha_0^2\pm 1) = 0
$;
\item[(ii)]
$
\ov{d}\equiv\mymod{b}{\mathfrak p}
$
since $\nu_{\mathfrak p}(\alpha_0\lambda^2 - 1)>0$; and
\item[(iii)]
$
\ov{d}\not= b
$
for all but finitely many $\mathfrak p$.
\end{itemize}
From  (i) we see that
$
P''_{a, b, \ov{c}, \ov{d}}
$
is well-defined over $\F_{\mathfrak p}$, and hence over $K$.  Its
$Y$-coordinate is zero modulo $\mathfrak p$, by (ii), but is non-zero in $K$
for all but finitely many $\mathfrak p$, by (iii).  By (ii),  the
discriminant (\ref{disc}) of
$
E''_{a, b, \ov{c}, \ov{d}}
$
is congruent modulo $\mathfrak p$ to
$
-27(a-\ov{c})^2(a^2 + a\ov{c} + \ov{c}^2)^2
$.
We saw that
$
a-\ov{c}\not\equiv\mymod{0}{\mathfrak p}
$,
and since
$
\alpha_0\lambda^2 \equiv\mymod{1}{\mathfrak p}
$,
$$
a^2 + a\ov{c} + \ov{c}^2
=
a^2
\Bigl(
  1 + \Bigl(\frac{\alpha_0^2\lambda^4}{\alpha_1^4}\Bigr)
    + \Bigl(\frac{\alpha_0^2\lambda^4}{\alpha_1^4}\Bigr)^2
\Bigr)
\equiv
\mymod{\frac{a^2}{\alpha_1^8}
\frac{1 - \alpha_1^{12}}
     {1 - \alpha_1^{4}}}
{\mathfrak p},
$$
which by (\ref{forwhich}) is a $\mathfrak p$-adic unit.  Thus
$
E''_{a, b, \ov{c}, \ov{d}}
$
is non-singular over $\F_{\mathfrak p}$, and hence over $K$.  Our earlier
comment about the $Y$-coordinate of 
$
P''_{a, b, \ov{c}, \ov{d}}
$
means that this point does not have order $2$ in 
$
E''_{a, b, \ov{c}, \ov{d}}(K)
$
but reduces to a point of order $2$ modulo $\mathfrak p$.    Apply Merel's theorem as
before and we see that
$
P''_{a, b, \ov{c}, \ov{d}}
$
is non-torsion in 
$
E''_{a, b, \ov{c}, \ov{d}}(K)
$.
\end{proof}

\vsp

\section{Proof of the Theorems}
     \label{sec:proof}

\begin{proof}[Proof of Theorem \ref{thm:pair}]
If the two elliptic curves in theorem \ref{thm:pair} are isomorphic over $K$,
then the desired positive rank twists follow readily from, for example, the
elementary construction mentioned after the statement of theorem
\ref{thm:zero}.
So suppose they are not isomorphic over $K$, in which case
lemma \ref{lem:tor} furnishes Weierstrass models
$
E_{a, b}
$
and
$
E_{\ov{c}, \ov{d}}
$
for them, such that
$
E''_{a, b, \ov{c}, \ov{d}}
$
is non-singular and contains a non-torsion rational point
$
P''_{a, b, \ov{c}, \ov{d}}
$.
That means
$
P'_{a, b, \ov{c}, \ov{d}}
$
is a non-torsion point on
$
E'_{a, b, \ov{c}, \ov{d}}
$.
For every integer $k$, denote by $x_k$ and $y_k$ the $x$- and $y$-coordinate
of
$
kP'_{a, b, \ov{c}, \ov{d}}
$.
Then $(x_k, y_k)$ is a rational solution to
$
f_{a, b}(x) = f_{\ov{c}, \ov{d}}(y)
$
for every $k$.  Denote by
$
z_k\in K
$
this  common value
$
f_{a, b}(x_k) = f_{\ov{c}, \ov{d}}(y_k)
$.
It is zero for at most finitely many $k$, so for all sufficiently large $k$,
each elliptic curve
$
z_k y^2 = f_{a, b}(x)
$
and
$
z_k y^2 = f_{\ov{c}, \ov{d}}(y)
$
has a rational point with non-zero $y$-coordinate.  At most finitely 
many quadratic twists of any elliptic curve over $K$ has $K$-rational torsion
points of order $>2$,
so the two twisted curves above have positive Mordell-Weil rank over $K$.
By lemma \ref{lem:genus},
for any $\lambda\in K^\times$ there are at most finitely
many $k$ for which 
$
z_k 
$
is $\lambda$ times a perfect square-free in $K$, and the theorem follows.
\end{proof}

\vsp

\begin{proof}[Proof of Theorem \ref{thm:zero}]
To say that both curves have $j$-invariant zero is to say that $a=c=0$.  If
in addition
$b=d$, then the two curves are identical, in which case the theorem follows
readily from the elementary construction mentioned in the introduction.
So from now on assume that $b\not=d$, in which case
$
E''_{0, b, 0, d}
$
becomes
$
E''': Y^2 = X^3 - 27 (b-d)^2/4
$.
The point $P'$ becomes $[1:1:0]$, which under the {\tt MAPLE} transformation
is the point at infinity for
$
E'''
$.
In general $E'''$ has no other $K$-rational point, so our argument above
 fails to generate a single quadratic twist, let alone infinitely many, for
which both
$
E_{0, b}
$
and
$
E_{0, d}
$
have positive rank.  On the other hand, we claim that there are infinitely
many
$
\lambda\in K
$
whose image in $K^\times/{K^\times}^3$ are pairwise distinct, such that the
cubic twist of
$
E''_{0, b, 0, d}
$
by $\lambda$ has positive Mordell-Weil rank.
Since this cubic twist can be written as
$
Y^2 = X^3 - 27 (\lambda b-\lambda d)^2 /4
$,
that means
$
E''_{0, \lambda b, 0, \lambda d}
$
has positive rank.  We can now resume the argument in the last paragraph
of the proof of theorem \ref{thm:pair} for the pair of sextic twists
$
E_{0, \lambda b}
$
and
$
E_{0, \lambda d}
$
(which only requires the \textit{existence} and not the explicit description
of a non-torsion point on
$
E''_{0, \lambda b, 0, \lambda d}
$)
and theorem \ref{thm:zero} follows.

\vsp

It remains to verify the claim.  The cubic twist
$
Y^2 = X^3 - 27 (\lambda b-\lambda d)^2 /4
$
is equivalent to
$
U^3 - V^3 = 4\lambda(d-b)
$.
Write $\lambda$ as
$
16(d-b)^2 t
$
and we are reduced to show that 
$
E_t: U^3 - V^3 = t
$
has positive rank for infinitely many $t$ whose image in
$
K^\times/{K^\times}^3
$
are pairwise distinct.
Since the $E_t$ are cubic twists, they have trivial torsion for all but
finitely many such $t$, so it suffices to show that $E_t$ has a non-trivial
rational
point for infinitely many $t$ whose image in
$
K^\times/{K^\times}^3
$
are pairwise distinct.
We proceed inductively as follow.  Pick $\mathfrak p\nmid 3$ in $\spec\oo_K$.
Since
$
U^3 - V^3 = (U-V)(U^2 + UV + V^2)
$,
if we choose $u, v\in\oo_K$ such that
$
u\equiv\mymod{(p+1)}{p^2}
$
and
$
v\equiv\mymod{1}{p^2}
$,
then $\mathfrak p|| (u^3 - v^3)$.  Thus
$
u^3 - v^3
$
could serve as one of the desired $t$ values.  Now, suppose we have
constructed finitely many such values $t_1, \ldots, t_m$.  Repeat the
process with a new
$
\mathfrak p_{m+1}\nmid (t_1 \cdots t_m)
$
and we obtain a new $t_{m+1}$; the conditions
$
\mathfrak p_{m+1} || t_{m+1}
$
and
$
\mathfrak p_{m+1}\nmid (t_1 \cdots t_m)
$
mean that the cube-free part of $t_{m+1}$ is different from those of the other
$t_i$.
Continue this process and we are done.
\end{proof}

\vsp

\begin{remark}
    \label{rem:triple}
If we try to treat the remaining case of theorem \ref{thm:pair} directly,
we need to
show that for 
infinitely many square-free integers $D$, the system
$$
(\gamma + \delta\sqrt{D})^3 + b - (\alpha + \beta\sqrt{D})^2 
\:=\:
0
\:=\:
(\gamma' + \delta'\sqrt{D})^3 + d - (\alpha' + \beta'\sqrt{D})^2 
$$
has a rational point 
$
(\alpha, \beta, \gamma, \delta, \alpha', \beta', \gamma', \delta')
$.
View each of the two equations above as the vanishing of an algebraic
integer in $\oo_K[\sqrt{D}]$, eliminate the variables $\alpha, \alpha'$ and $D$
and we arrive at a family of curves of geometric genus $11$ in the variables
$
\gamma, \gamma'
$
over the affine parameters
$
\beta, \beta', \delta, \delta'
$.
Showing that the total space of such a family has one rational point,
let alone infinitely many, with $bd\not=0$ seems highly non-trivial.

\vsp

If $K$ contains a primitive third root of unity $\zeta_3$, then the projective
curve
$
E_{a, b, c, d}
$
acquires the additional $K$-rational points
$
[1: \zeta_3: 0]
$
and
$
[1: \zeta_3^2: 0]
$.
However, when $a=c=0$ we check that these become $2$-torsion points on
$
E_{a, b, c, d}(K)
$,
and hence they cannot be use to generate an infinite collection of
positive rank twists over $K$.
\end{remark}

\vsp

\begin{ack}
Part of this work was carried out during the Banff workshop on analytic methods
on diophantine equations.  We would like to thank to the organizers for
providing a conducive working environment, and Professors Hajir, Rohrlich
and Rosen
for 
useful discussion.
This research is supported in part by NSA grant H98230-05-1-0069.
\end{ack}

\vtsp

\bibliographystyle{amsalpha}

\vfill
\hrule
\end{document}